\documentclass[a4paper,12pt]{article}

\usepackage[text={6.5in,8in},centering]{geometry}

\usepackage{graphicx,url,psfrag,subfig}
\RequirePackage[OT1]{fontenc}
\RequirePackage{amsthm,amsmath,amssymb,amscd}
\RequirePackage[colorlinks,citecolor=blue,urlcolor=blue]{hyperref}
\usepackage{nomencl}
\usepackage{enumerate}
\usepackage{multirow}
\usepackage{tensor}
\usepackage{etex}

\usepackage{caption}

\usepackage{svn-multi}[2009/03/03]

\theoremstyle{plain}
\newtheorem{theorem}{Theorem}[section]
\newtheorem{corollary}[theorem]{Corollary}

\newtheorem{proposition}[theorem]{Proposition}

\theoremstyle{definition}
\newtheorem{definition}[theorem]{Definition}

\newtheorem{remark}[theorem]{Remark}

\newtheorem{example}[theorem]{Example}

\numberwithin{equation}{section}

\allowdisplaybreaks

\svnidlong
{$HeadURL: file:///home/lechen/svn/Article-Heat/HeatAnnounce.tex $}
{$LastChangedDate: 2012-10-04 17:49:18 +0200 (Thu, 04 Oct 2012) $}
{$LastChangedRevision: 106 $}
{$LastChangedBy: lechen $}
\svnRegisterAuthor{Le}{Le Chen}

\newcommand{\E}{\mathbb{E}}
\newcommand{\W}{\dot{W}}
\newcommand{\ud}{\ensuremath{\mathrm{d}}}
\newcommand{\Ceil}[1]{\lceil #1 \rceil}

\newcommand{\Norm}[1]{\left|\left|  #1   \right|\right|}

\newcommand{\calB}{\mathcal{B}}

\newcommand{\calF}{\mathcal{F}}

\newcommand{\calK}{\mathcal{K}}
\newcommand{\calH}{\mathcal{H}}

\newcommand{\calM}{\mathcal{M}}
\newcommand{\calN}{\mathcal{N}}

\newcommand{\bbN}{\mathbb{N}}

\newcommand{\R}{\mathbb{R}}

\newcommand{\Erf}{\ensuremath{\mathrm{erf}}}
\newcommand{\Erfc}{\ensuremath{\mathrm{erfc}}}

\DeclareMathOperator{\Lip}{\mathit{L}}
\DeclareMathOperator{\LIP}{Lip}
\DeclareMathOperator{\lip}{\mathit{l}}

\DeclareMathOperator{\Vip}{\overline{\varsigma}}
\DeclareMathOperator{\vip}{\underline{\varsigma}}
\DeclareMathOperator{\vv}{\varsigma}
\DeclareMathOperator{\sd}{\beta}

\title{{\bf The nonlinear stochastic heat equation \\
with rough initial data: \\
a summary of some new results}}

\author{{\bf Le Chen$^*$} and {\bf Robert C. Dalang\footnote{Research
partially supported by the Swiss National Foundation for Scientific
Research.}}\\
\\
\it\small
Institut de math\'ematiques\\
\it\small \'Ecole Polytechnique F\'ed\'erale\\
\it\small Station 8 \\
\it\small CH-1015 Lausanne\\
\it\small Switzerland\\
\small \textit{e-mails:}
le.chen@epfl.ch, robert.dalang@epfl.ch
}
\date{}

\begin{document}
\maketitle
\begin{center}
\begin{minipage}[rct]{5 in}
\footnotesize
\textbf{Abstract:} This is a preliminary announcement of results in the PhD.
thesis
of the first author concerning the nonlinear stochastic heat equation in the
spatial
domain $\R$, driven by space-time white noise. A central special case
is the parabolic Anderson model. The initial condition is taken to be a measure
on $\R$, such as the Dirac delta function, but this measure may also have
non-compact support and even be non-tempered (for instance with exponentially
growing tails). Existence and uniqueness is proved without appealing to
Gronwall's lemma, by keeping tight control over moments in the Picard iteration
scheme.  Upper and lower bounds on all $p$-th moments $(p\ge 2)$ are obtained.
These bounds become equalities for the parabolic Anderson model when $p=2$.
The growth indices introduced by Conus and Khoshnevisan
\cite{ConusKhosh10Farthest} are determined and, despite the irregular initial
conditions, H\"older continuity of the solution for $t>0$ is established.

\vspace{2ex}
\textbf{AMS 2010 subject classifications:}
Primary 60H15. Secondary 60G60, 35R60.

\vspace{2ex}
\textbf{Keywords:}
nonlinear stochastic heat equation, parabolic Anderson model, rough initial
data, growth indices.
\vspace{4ex}
\end{minipage} 
\end{center}



\graphicspath{{figs/}}

\section{Introduction}
The stochastic heat equation
\begin{align}\label{E2:Heat}
\begin{cases}
\displaystyle \left(\frac{\partial }{\partial t} - \frac{\nu}{2}
\frac{\partial^2 }{\partial x^2}\right) u(t,x) =  \rho(u(t,x))
\:\W(t,x),&
x\in \R,\; t \in\R_+^*, \\
\displaystyle \quad u(0,\cdot) = \mu(\cdot)\;,
\end{cases}
\end{align}
\nomenclature{$\W$}{Space-time white noise \nomrefpage}\noindent
where $\W$ is space-time white noise, $\rho(u)$ is globally Lipschitz,
$\mu$ is the initial data, and $\R_+^*=\;]0,+\infty[$, has been intensively
studied during last two decades by many authors
\cite{BertiniCancrini94Intermittence,CarmonaMolchanov94,ConusEct12Initial,
ConusKhosh10Farthest,
ConusKhosh10Weak,DalangKhNualart07HittingAdditive,DalangKhNualart09HittingMult,
FoondunKhoshnevisan08Intermittence,Mueller91Support,PospisilTribe07Parameter,
SanSoleSarra99Holder, SanSoleSarra99Path,Shiga94Two,Walsh86} .
In particular, the special
case $\rho(u)=\lambda u$ is called {\it the parabolic Anderson model}
\cite{CarmonaMolchanov94}. Our
work focuses on this equation with general deterministic initial data, and we
study
how the initial data affects the solution.

The one-dimensional heat kernel function is
\begin{align}\label{E2:G-Heat}
G_\nu(t,x) := \frac{1}{\sqrt{2\pi \nu t}} \exp\left\{-\frac{|x|^2}{2\nu
t}\right\},\quad (t,x)\in\R_+^*\times\R\:.
\end{align}
\nomenclature{$G_\nu(t,x)$}{heat kernel function \nomrefeqpage}\noindent
For the existence of random field solutions to \eqref{E2:Heat}, the case where
the initial data
$\mu$ is a bounded and measurable function is covered by the classical theory of
Walsh \cite{Walsh86}.
When $\mu$ is a positive Borel measure on $\R$ such that
\begin{align}\label{E2:BC-InitD}
\sup_{t\in [0,T]} \sup_{x\in\R} \sqrt{t}
\left(\mu*G_\nu(t,\circ)\right)(x)<\infty,\quad\text{for all $T>0$},
\end{align}
where $*$ denotes convolution in the spatial variable, Bertini and Cancrini
\cite{BertiniCancrini94Intermittence} gave an ad-hoc definition for the
Anderson model via a smoothing of the space-time white
noise and a Feynman-Kac type formula. Their analysis depends heavily on
properties of the local times of Brownian bridges. Recently, Conus and
Khoshnevisan \cite{ConusKhosh10Weak} constructed a weak solution defined through
certain norms on random fields. The initial data has to verify certain
technical conditions, which include the Dirac delta function in some of
their cases.
In particular, the solution is defined for almost all $(t,x)$, but not at
specific $(t,x)$. More recently, Conus, Joseph, Khoshnevisan and Shiu
\cite{ConusEct12Initial} also studied random field
solutions. In particular, they require the initial data to be a finite measure
of compact support. We improve the existence result by working under a
much weaker condition on initial data, namely, $\mu$ can be any signed
Borel measure over $\R$ such that
\begin{align}\label{E2:J0finite}
\left(|\mu| * G_\nu(t,\cdot)\right)
(x)<+\infty\;, \quad \text{for all $t>0$ and $x\in\R$}\;,
\end{align}
\nomenclature{$*$}{space convolution \nomrefeqpage}\noindent
where, from the Jordan decomposition, $\mu=\mu_+-\mu_-$ where
$\mu_\pm$ are two non-negative Borel measures
and $|\mu|:= \mu_++\mu_-$.
For instance, if $\mu(\ud x)=f(x)\ud x$, then $f(x)=\exp\left(a|x|^p\right)$,
$a>0$, $p\in \;]0,2[$, (i.e.,
exponential growth at $\pm \infty$), will satisfy this condition.
Proposition \ref{P2:D-Delta} below shows that
initial data cannot be extended beyond measures to other Schwartz distributions,
even with compact support.

Moreover, we obtain estimates for the moments $\E(|u(t,x)|^p)$ with both $t$ and
$x$ fixed for all even integers $p\ge 2$. In particular, for the
parabolic Anderson model, we give an explicit formula for the second moment of
the solution.
When the initial data is either Lebesgue measure or the Dirac delta function,
we
give explicit formulae for the two-point correlation functions (see
\eqref{E2:TP-Lebesgue} and \eqref{E2:TP-Delta} below), which
can be compared to the integral form in Bertini and Cancrini's paper
\cite[Corollaries 2.4 and
2.5]{BertiniCancrini94Intermittence}.

Our proof of existence is based on the standard Picard iteration. The main
difference from the conventional situation is that instead of applying
Gronwall's lemma to bound the second moment from above, we show that the
sequence of the second moments in the Picard iteration converges to an
explicit formula (in the case of the parabolic Anderson model).

After establishing the existence of random field solutions, we study whether
the solution
exhibits intermittency properties. More precisely, define the {\it upper and
lower Lyapunov exponents} for constant initial data (Lebesgue measure) as
follows
\begin{align}\label{E2:Lyapunov}
 \overline{\lambda}_p:=& \mathop{\lim\sup}_{t\rightarrow+\infty}\frac{\log
\E\left[|u(t,x)|^p\right]}{t},\qquad
\underline{\lambda}_p:=\mathop{\lim\inf}_{t\rightarrow+\infty}\frac{\log
\E\left[|u(t,x)|^p\right]}{t} \;.
\end{align}
\nomenclature{$\overline{\lambda}_p$}{upper Lyapunov exponent of
order $p$ \nomrefeqpage}\noindent
\nomenclature{$\underline{\lambda}_p$}{lower Lyapunov exponent of
order $p$ \nomrefeqpage}\noindent
Following Bertini and Cancrini \cite{BertiniCancrini94Intermittence}, we
say that the solution is {\it intermittent} if
$\lambda_n:=\underline{\lambda}_n =\overline{\lambda}_n$ and the strict
inequalities
\begin{align}\label{E2:Intermit}
\lambda_1 <
\frac{\lambda_2}{2}<\cdots<\frac{\lambda_n}{n}<\cdots
\end{align}
are satisfied. Carmona and Molchanov gave the following definition
\cite[Definition III.1.1, on p. 55]{CarmonaMolchanov94}:

\begin{definition}[Intermittency]\label{D2:Intermit}
Let $p$ be the smallest integer for which $\lambda_p>0$. When $p<\infty$, we
say that the solution $u(t,x)$ shows {\em (asymptotic) intermittency of
order $p$} and {\em full intermittency} when $p=2$.
\end{definition}

They showed that full intermittency implies the
intermittency defined by
\eqref{E2:Intermit} (see \cite[III.1.2, on p. 55]{CarmonaMolchanov94}).
This mathematical definition of intermittency is related to the property that
the solutions develop
high peaks on some small ``islands".
The parabolic Anderson model has been well
studied:
see \cite{CarmonaMolchanov94} for a discrete approximation and
\cite{BertiniCancrini94Intermittence} and
\cite{FoondunKhoshnevisan08Intermittence} for the continuous version.
Further discussion can be found in \cite{Zeldovich90Almighty}.

When the initial data are not homogeneous, in particular, when they have certain
exponential decrease at infinity, Conus and Khoshnevisan
\cite{ConusKhosh10Farthest} defined the following
{\it lower and upper exponential growth indices}:
\begin{align}
\label{E2:GrowInd-0}
\underline{\lambda}(n):= &
\sup\left\{\alpha>0: \underset{t\rightarrow \infty}{\lim\sup}
\frac{1}{t}\sup_{|x|\ge \alpha t} \log \E\left(|u(t,x)|^n\right) >0
\right\}\;,\\
\label{E2:GrowInd-1}
 \overline{\lambda}(n) := &
\inf\left\{\alpha>0: \underset{t\rightarrow \infty}{\lim\sup}
\frac{1}{t}\sup_{|x|\ge \alpha t} \log \E\left(|u(t,x)|^n\right) <0
\right\}\:,
\end{align}
\nomenclature{$\underline{\lambda}(n)$}{lower exponential growth indices of
order $n$ \nomrefeqpage}\noindent
\nomenclature{$\overline{\lambda}(n)$}{upper exponential growth indices of
order $n$ \nomrefeqpage}\noindent
and proved that if the initial data $\mu$ is a non-negative, lower
semicontinuous function with compact support of positive measure, then
for the Anderson model ($\rho(u)=\lambda u$),
\[
\frac{\lambda^2}{2\pi} \le \underline{\lambda}(2)
\le \overline{\lambda}(2) \le \frac{\lambda^2}{2}\;.
\]
We improve this result by
showing that $\underline{\lambda}(2)= \overline{\lambda}(2) =
\lambda^2/2$, and extend this to more general measure-valued initial data.
This is possible mainly thanks to our explicit formula for the second moment.

We now discuss the regularity of the random field solution.
Denote by $C_{\beta_1,\beta_2}(D)$ the set of random fields whose
trajectories are almost surely $\beta_1$-H\"older continuous in time and
$\beta_2$-H\"older continuous in space on the domain $D\subseteq
\R_+\times\R$, and let
\[
C_{\beta_1-,\beta_2-}(D) := \bigcup_{\alpha_1\in \;\left]0,\beta_1\right[}
\bigcup_{\alpha_2\in \;\left]0,\beta_2\right[} C_{\alpha_1,\alpha_2}(D)\;.
\]
\nomenclature{$C_{\beta_1,\beta_2}(D)$}{ the set of random fields whose
trajectories are almost surely $\beta_1$-H\"older continuous in time and
$\beta_2$-H\"older continuous in space over domain $(t,x)\in D\subseteq
\R_+\times\R$ \nomrefpage}\noindent
In Walsh's notes \cite[Corollary 3.4, p.
318]{Walsh86}, a slightly different equation was studied and the H\"older
exponents given (for both space and time) are $1/4- \epsilon$.
Bertini and Cancrini \cite{BertiniCancrini94Intermittence} stated in their
paper that the random field solution for the parabolic Anderson model with
initial data satisfying \eqref{E2:BC-InitD} belongs to
$C_{\frac{1}{4}-,\frac{1}{2}-}(\R_+^*\times\R)$.
In \cite{PospisilTribe07Parameter,Shiga94Two}, the authors showed that if the
initial data is a continuous function with certain exponentially growing tails,
then
\begin{align}\label{E2:Shiga}
u\in C_{\frac{1}{4}-,\frac{1}{2}-}(\R_+\times\R),\quad\text{a.s.}
\end{align}
Sanz-Sol\'e and Sarr\`a \cite{SanSoleSarra99Holder} considered the stochastic
heat equation over $\R^d$ with spatially homogeneous colored noise which is
white in time. Let $\tilde{\mu}$ be the spectral measure satisfying
\[
\int_{\R^d}\frac{\tilde{\mu}(\ud
\xi)}{\left(1+|\xi|^2\right)^\eta}<+\infty,
\quad\text{for some $\eta\in \;]0,1[$.}
\]
They proved that if the initial data is a bounded
$\rho$-H\"older continuous function for some $\rho\in \;]0,1[$, then the
solution
is in
\[
u(t,x)\in C_{\frac{1}{2}(\rho\wedge (1-\eta))-, \rho\wedge
(1-\eta)-}
\left(\R_+^*\times\R\right)\;.
\]
For the case of space-time white noise, the spectral measure $\tilde{\mu}$
is Lebesgue measure and hence $\eta$ can be $1/2-\epsilon$ for any
$\epsilon>0$. Their result (\cite[Theorem 4.3]{SanSoleSarra99Path})
reduces to
\[
u(t,x)\in C_{\left(\frac{1}{4}\wedge\frac{\rho}{2}\right)-,
\left(\frac{1}{2}\wedge\rho\right)-}\left(\R_+^*\times\R\right)\;.
\]
More recently, Conus {\it et al} proved in their paper
\cite[Lemma 9.3]{ConusEct12Initial} that the random field solution is H\"older
continuous in $x$ with exponent $1/2-\epsilon$ (for initial data that is a
finite measure). They did not give the regularity estimate over the time
variable. In their papers \cite{DalangKhNualart07HittingAdditive,
DalangKhNualart09HittingMult}, 
Dalang, Khoshnevisan and Nualart considered a system of heat equations with
vanishing initial conditions subject to space-time white noise, and proved that
the solution is jointly H\"older continuous with exponents $1/4-$ in time and
$1/2-$ in space.
We extend the $C_{\frac{1}{4}-,\frac{1}{
2}-}\left(\R_+^*\times\R\right)$-H\"older
continuity result to measure-valued initial data satisfying \eqref{E2:J0finite}.
We show that the result in \eqref{E2:Shiga} should exclude the time line
$t=0$ unless the initial data $\mu$ is $1/2$-H\"older continuous.

The difficulties for the proof of the H\"older continuity of the random field
solution lie in the fact that for the initial data satisfying
\eqref{E2:J0finite}, the $p$-th moment $\E\left[|u(t,x)|^p\right]$ is neither
bounded for $x\in\R$, nor for $t\in [0,T]$. Standard
techniques, which isolate the effects of initial data by the
$L^p(\Omega)$-boundedness of the solution, fail in our case. Instead, the
initial data play an active role in our proof.
Note that Fourier transforms are not applicable here because $\mu$ need not be
a tempered measure.

\section{Main Results}\label{S2:MainRes}

Denote the solution to the homogeneous equation
\begin{align}\label{E2:Heat-home}
\begin{cases}
\displaystyle \left(\frac{\partial }{\partial t} - \frac{\nu}{2}
\frac{\partial^2 }{\partial x^2}\right) u(t,x) =  0,&
x\in \R,\; t \in\R_+^*, \\
\displaystyle \quad u(0,\cdot) = \mu(\cdot)\;,
\end{cases}
\end{align}
by
\[
J_0(t,x) := \left(\mu* G_{\nu}(t,\cdot)\right)(x) = \int_\R G_\nu(t,x-y)\mu(\ud
y)\:,\quad (t,x)\in\R_+^*\times\R\:.
\]
\nomenclature{$J_0(t,x)$}{solutions to the homogeneous equation
\nomrefpage}\noindent
Note that $J_0(t,x)$ is well-defined by the hypothesis \eqref{E2:J0finite}.
We formally  rewrite the stochastic partial differential
equation \eqref{E2:Heat} in the integral form (mild form):
\begin{align}\label{E2:WalshSI}
 u(t,x) &= J_0(t,x) + I(t,x)
\end{align}
where
\begin{align}\label{E2:WalshSI-I}
I(t,x):=\iint_{[0,t]\times\R} G_\nu(t-s,x-y)
\rho\left( u(s,y) \right)
\W(\ud s,\ud y)\:.
\end{align}
By convention, $I(0,x)=0$.
The above stochastic integral
is defined in the sense of Walsh \cite{Walsh86,DalangEtc08Minicourse}.

\subsection{Notations and Conventions}
Assume that the function $\rho:\R\mapsto \R$ is globally Lipschitz
continuous with Lipschitz constant $\LIP_\rho>0$.
\nomenclature{$\LIP_\rho$}{Lipschitz constant \nomrefpage}
We need some growth conditions on $\rho$:
Assume that
\begin{align}\label{E2:LinGrow}
|\rho(x)|^2 \le \Lip_\rho^2 \left(\Vip^2 +|x|^2\right),\qquad \text{for
all $x\in\R$}\;,
\end{align}
for some constants $\Lip_\rho>0$ and $\Vip \ge 0$.
\nomenclature{$\Lip_\rho$}{linear growth constant (upper bound)
\nomrefeqpage}
\nomenclature{$\Vip$}{constant related to the linear growth (upper
bound) \nomrefeqpage}
Note that $\sqrt{2}\LIP_\rho\le \Lip_\rho$, and the inequality may be strict.
In order to bound the second moment from below,
we will sometimes assume that for some constants $\lip_\rho>0$ and $\vip \ge 0$,
\begin{align}\label{E2:lingrow}
|\rho(x)|^2\ge \lip_{\rho}^2\left(\vip^2+|x|^2\right),\qquad \text{for
all $x\in\R$}\;.
\end{align}
\nomenclature{$\lip_\rho$}{linear growth constant (lower bound)
\nomrefeqpage}\noindent
\nomenclature{$\vip$}{constant related to the linear growth (lower bound)
\nomrefeqpage}\noindent
We shall also give special attention to the linear case (the parabolic
Anderson model) $\rho(u)=\lambda u$ with $\lambda\ne 0$, which is a
special case of the following quasi-linear growth condition:
for some constant $\vv\ge 0$,
\begin{align}\label{E2:qlinear}
|\rho(x)|^2= \lambda^2\left(\vv^2+|x|^2\right),\qquad  \text{for
all $x\in\R$}\;.
\end{align}
\nomenclature{$\vv$}{constant related to the parabolic Anderson model
\nomrefeqpage}

We use the convention that $G_\nu(t,\cdot)\equiv 0$ if $t\le 0$. Hence,
the integral region in the stochastic integral in \eqref{E2:WalshSI-I} can be
written as $\R_+\times\R$.

Define a kernel function
\begin{align} \label{E2:K}
 \calK\left(t,x;\nu,\lambda\right) := G_{\frac{\nu}{2}}(t,x)
\left(\frac{\lambda^2}{\sqrt{4\pi\nu t}}+\frac{\lambda^4}{2\nu}
\: e^{\frac{\lambda^4 t}{4\nu}}\Phi\left(\lambda^2
\sqrt{\frac{t}{2\nu}}\right)\right) \;,
\end{align}
\nomenclature{$\calK\left(t,x;\nu,\lambda\right)$}{a kernel function
\nomrefeqpage}\noindent
for all $(t,x)\in\R_+^*\times \R$,
where $\Phi(x)$ is the probability distribution function of the standard normal
distribution:
\[
\Phi(x) := \int_{-\infty}^x \frac{e^{-y^2/2}}{\sqrt{2\pi}} \ud y\;.
\]
\nomenclature{$\Phi(x)$}{the probability distribution function of the standard
normal distribution \nomrefpage}\noindent
We also use the error function
$\Erf(x):=\frac{2}{\sqrt{\pi}}\int_0^x e^{-y^2}\ud y$ and its complement
$\Erfc(x):=1-\Erf(x)$.
Clearly,
\begin{align*}
\Phi(x) =
\frac{1}{2}\left(1+\Erf\left(x/\sqrt{2}\right)\right)\;,
\end{align*}
\begin{align*}
\Erf(x) &= 2\Phi\left(\sqrt{2}\:x\right)-1,\quad
\Erfc(x) = 2\left(1-\Phi\left(\sqrt{2}\:x\right)\right)\;.
\end{align*}
\nomenclature{$\Erf(x)$}{ the error function
$\Erf(x):=\frac{2}{\sqrt{\pi}}\int_0^x e^{-x^2}\ud x$.\nomrefpage}
\nomenclature{$\Erfc(x)$}{ the complementary error function $\Erfc(x) =
1-\Erf(x)$ \nomrefpage}

We use $\star$ to denote the simultaneous convolution in both space
and time variables.
\nomenclature{$\star$}{space-time convolution \nomrefpage}
Define another function
\begin{align}\label{E2:H}
\calH(t;\nu,\lambda):=\left(1\star \calK \right)(t,x) = 2 e^{\frac{\lambda^4\:
t}{4\nu}}\Phi\left(\lambda^2\sqrt{\frac{t}{2\nu}}\right)-1\;.
\end{align}
Clearly, $\calK\left(t,x;\nu,\lambda\right)$ can be written as
\[
\calK\left(t,x;\nu,\lambda\right)=G_{\nu/2}(t,x)
\left(
\frac{\lambda^2}{\sqrt{4\pi\nu t}}+ \frac{\lambda^4}{4\nu}
\left[\calH(t;\nu,\lambda)+1\right]
\right)
\;.
\]
We use the following conventions:
\begin{align}\label{E2:K-cm}
\calK(t,x) &:=  \calK\left(t,x \:
;\: \nu,\lambda\right)\;, \\
\label{E2:upperK}
\overline{\calK}(t,x)   &:=  \calK\left(t,x \:
;\: \nu,\Lip_\rho\right)\;,\\
\label{E2:lowerK}
\underline{\calK}(t,x) &:=\calK\left(t,x \:;\: \nu,\lip_\rho \right)\;,\\
\label{E2:hatK}
\widehat{\calK}_p(t,x) &:=\calK\left(t,x \:
;\: \nu,a_{p,\Vip}\:z_p\: \Lip_\rho \right)\;,\quad\text{for all $p>2$}\;,
\end{align}
where $z_p$ is the universal constant in the
Burkholder-Davis-Gundy inequality (in particular, $z_2=1$) 
and $a_{p,\Vip}$ is a constant defined as
\begin{align}\label{E2:a_pv}
a_{p,\Vip} \: :=\:
\begin{cases}
2^{(p-1)/p}& \text{if $\Vip\ne 0,\: p>2$,}\cr
\sqrt{2} & \text{if $\Vip =0,\: p>2$,}\cr
1 & \text{if $p=2$}.
\end{cases}
\end{align}
We only need  to keep in mind that
$a_{p,\Vip}\le 2$.
Note that the kernel function $\widehat{\calK}_p(t,x)$ implicitly depends on
$\Vip$ through $a_{p,\Vip}$ which will be clear from the context. If $p=2$, then
$\widehat{\calK}_2(t,x) = \overline{\calK}(t,x)$.

Similarly $\overline{\calH}(t)$,
$\underline{\calH}(t)$ and
$\widehat{\calH}_p(t)$ denote the kernel functions
with $\lambda$ in $\calH(t)$ replaced by
$\Lip_\rho$, $\lip_\rho$ and $a_{p,\Vip} z_p \Lip_\rho$,
respectively. Again $\widehat{\calH}_p(t)$ depends on $\Vip$ implicitly which
will be clear from the context.

Let us set up the filtered probability space. Let
\[
\Big\{\: W_t(A):\:A\in\calB_b\left(\R\right),t\ge 0 \: \Big\}
\]
be a space-time white noise
defined on a probability space $(\Omega,\calF,P)$, where
$\calB_b\left(\R\right)$ is the
collection of Borel measurable sets with finite Lebesgue measure.
Let  $(\calF_t,t\ge 0)$ be the standard filtration generated by this space-time
white noise. More precisely, let
\[
\calF_t^0 := \sigma\left(W_s(A):\: 0\le s\le t,
A\in\calB_b\left(\R\right)\right)\vee
\calN,\quad t\ge 0
\]
be the natural filtration augmented by all $P$-null sets in $\calF$.
Define $\calF_t := \calF_{t+}^0 = \wedge_{s>t}\calF_s^0$ for any $t\ge 0$.
In the following, we fix this filtered
probability space $\left\{\Omega,\calF,\{\calF_t:t\ge0\},P\right\}$.
We use $\Norm{\cdot}_p$ to denote the
$L^p(\Omega)$-norm.
\nomenclature{$\Norm{\cdot}_p$}{ the $L^p(\Omega)$-norm \nomrefpage}
Denote $\Ceil{p}_2:=2\Ceil{p/2}$, which is the
smallest even integer greater
than or equal to $p$.
\nomenclature{$\Ceil{p}_2$}{the smallest even integer greater than or equal to
$p$ \nomrefpage}

Let $\calM(\R)$ be the set of locally finite Borel measures over $\R$.
\nomenclature{$\calM(\R)$}{the set of locally finite Borel measures on $\R$
\nomrefpage}
Define
\begin{align}\label{E2:MGa}
\calM_G^{\sd}(\R):= \left\{\mu\in\calM(\R): \: \int_\R e^{\sd|x|}|\mu|
(\ud x)<+\infty\right\}\;, \quad \sd\ge 0,
\end{align}
\nomenclature{$\calM_G^{\sd}(\R)$}{the set of Borel measures on $\R$
that have exponential decays at infinite \nomrefeqpage}\noindent
where $|\mu|=\mu_++\mu_-$ is the Jordan decomposition of a measure into two
non-negative measures.
We use subscript ``$+$'' to denote the subset of non-negative measures.
For example, $\calM_+(\R)$ is the set of non-negative Borel measures over $\R$
and
$\calM_{G,+}^{\sd}(\R) = \calM_G^{\sd}(\R) \cap \calM_+(\R)$.

A random field $Y=\left(Y(t,x):\: (t,x)\in \R_+^*\times
\R\right)$ is said to be {\it $L^p(\Omega)$-continuous}, $p\ge 2$, if
for all $(t,x) \in \R_+^*\times\R$, 
\[
\lim_{\left(t',x'\right)\rightarrow (t,x)}\Norm{Y(t,x)-Y\left(t',x'\right)}_p
=0\;.
\]

\subsection{Existence, Uniqueness and Moments}
We first give the definition of the random field solution as follows:
\begin{definition}\label{D2:Solution}
A solution $u=\left(u(t,x):(t,x)\in\R_+^*\times\R\right)$ to \eqref{E2:Heat} (or
\eqref{E2:WalshSI}) is called a {\it
random field solution} if
\begin{enumerate}[(1)]
 \item $u$ is adapted, i.e., for all
$(t,x)\in\R_+^*\times\R$, $u(t,x)$ is
$\calF_t$-measurable;
\item $u$ is jointly measurable with respect to
$\calB\left(\R_+^*\times\R\right)\times\calF$;
\item $\left(G_\nu^2 \star \Norm{\rho(u)}_2^2\right)(t,x)<+\infty$
for all $(t,x)\in\R_+^*\times\R$, and
the function
$(t,x)\mapsto I(t,x)$ mapping from $\R_+^*\times\R$ into
$L^2(\Omega)$ is continuous;
\item $u$ satisfies \eqref{E2:Heat} (or \eqref{E2:WalshSI}) almost surely,
for all $(t,x)\in\R_+^*\times\R$.
\end{enumerate}
\end{definition}

The first main result is stated as follows.
\begin{theorem}[Existence, uniqueness, moments]
\label{T2:ExUni}
Suppose that
\begin{enumerate}[(i)]
 \item the initial data $\mu$ is a signed Borel measure such that
\eqref{E2:J0finite} holds;
 \item the function $\rho$ is Lipschitz continuous such that the linear growth
condition \eqref{E2:LinGrow} holds.
\end{enumerate}
Then the stochastic integral equation \eqref{E2:WalshSI} has a random field
solution $u=\{u(t,x): t>0,x\in\R\}$ (note that $t>0$) in the sense of
Definition \ref{D2:Solution}. This solution has the
following properties:
\begin{enumerate}[(1)]
\item $u$ is unique (in the sense of versions);
\item $(t,x)\mapsto u(t,x)$ is $L^p(\Omega)$-continuous for all integers $p\ge
2$;
\item For all even integers $p\ge 2$, the $p$-th moment of the solution
$u(t,x)$ satisfies the upper bound
\begin{align}
\label{E2:SecMom-Up}
\Norm{u(t,x)}_p^2 \le
\begin{cases}
J_0^2(t,x) + \left(J_0^2\star \overline{\calK} \right) (t,x) +
\Vip^2 \: \overline{\calH}(t),& \text{if $p=2$,}\cr \cr
2J_0^2(t,x) + \left(2J_0^2\star \widehat{\calK}_p \right) (t,x) +
\Vip^2 \: \widehat{\calH}_p(t),& \text{if $p>2$,}
\end{cases}
\end{align}
for all $t>0$, $x\in\R$,
and the two-point correlation satisfies the upper bound
\begin{multline}
\label{E2:TP-Up}
\E\left[u(t,x)u\left(t,y\right)\right] \\ \le
J_0(t,x)J_0\left(t,y\right)+\Lip_\rho^2\int_0^t\ud s \int_{\R}
\overline{f}(s,z) G_\nu(t-s,x-z) G_\nu(t-s,y-z) \ud z\\
+ \frac{\Lip_\rho^2\Vip^2}{\nu}|x-y| \left(\Phi\left(\frac{|x-y|}{\sqrt{2\nu
t}}\right)-1\right) +
2 \Lip_\rho^2 \Vip^2 t\: G_{2\nu}(t,x-y)\;,
\end{multline}
for all $t>0$, $x,y\in\R$, where $\overline{f}(s,z)$ denotes the right
hand side of \eqref{E2:SecMom-Up} for $p=2$;
\item If $\rho$ satisfies \eqref{E2:lingrow}, then the
second moment satisfies the lower bound
\begin{align}
\label{E2:SecMom-Lower}
\Norm{u(t,x)}_2^2 \ge J_0^2(t,x) + \left(J_0^2 \star \underline{\calK} \right)
(t,x)+
\vip^2\: \underline{\calH}(t)
\end{align}
for all $t>0$, $x\in\R$, and the two-point correlation satisfies the lower bound
\begin{multline}
\label{E2:TP-Lower}
\E\left[u(t,x)u\left(t,y\right)\right] \\ \ge
J_0(t,x)J_0\left(t,y\right) +\lip_\rho^2\int_0^t\ud s \int_{\R}
\underline{f}(s,z)
G_\nu(t-s,x-z) G_\nu(t-s,y-z) \ud z\\
+ \frac{\lip_\rho^2 \vip^2}{\nu}|x-y| \left(\Phi\left(\frac{|x-y|}{\sqrt{2\nu
t}}\right)-1\right) +
2 \lip_\rho^2 \vip^2 t\: G_{2\nu}(t,x-y)\:,
\end{multline}
for all $t>0$, $x,y\in\R$, where $\underline{f}(s,z)$ denotes the right
hand side of \eqref{E2:SecMom-Lower};
\item In particular, for the quasi-linear case $|\rho(u)|^2=\lambda^2
\left(\vv^2+u^2\right)$, the second moment has the explicit expression
\begin{align}
 \label{E2:SecMom}
 \Norm{u(t,x)}_2^2 = J_0^2(t,x) + \left(J_0^2 \star \calK \right) (t,x)+ \vv^2
\: \calH(t)\;,
\end{align}
for all $t>0$, $x\in\R$,  and the two-point correlation is given by
\begin{multline}
\label{E2:TP}
\E\left[u(t,x)u\left(t,y\right)\right] \\ =
J_0(t,x)J_0\left(t,y\right) +\lambda^2 \int_0^t\ud s \int_{\R} f(s,z)
G_\nu(t-s,x-z) G_\nu(t-s,y-z) \ud z\\
+ \frac{\lambda^2 \vv^2}{\nu}|x-y| \left(\Phi\left(\frac{|x-y|}{\sqrt{2\nu
t}}\right)-1\right) +
2 \lambda^2 \vv^2 t\: G_{2\nu}(t,x-y)\:,
\end{multline}
for all $t>0$, $x,y\in\R$, where $f(s,z)=\Norm{u(s,z)}_2^2$ is defined in
\eqref{E2:SecMom}.
\end{enumerate}
\end{theorem}

\begin{corollary}[Lebesgue initial data] \label{C2:TP-Lebesgue}
Suppose that $|\rho(u)|^2=\lambda^2(\vv^2+u^2)$ and $\mu$ is
Lebesgue measure. Then for all $t>0$ and $x,y\in\R$,
\begin{multline}
\label{E2:TP-Lebesgue}
\E\left[u(t,x)u\left(t,y\right)\right] 
=
1+(1+\vv^2)
\Bigg(\exp\left(\frac{\lambda ^4 t-2 \lambda ^2 |x-y|}{4
   \nu }\right) \\
\times \Erfc\left(\frac{|x-y|-\lambda ^2 t}{2
   \sqrt{\nu  t}}\right)
-
\Erfc\left(\frac{|x-y|}{2 \sqrt{\nu t}}\right)
\Bigg)\;.
\end{multline}
In particular, when $y=x$, we have
\begin{align}
\label{E2:SecMom-Lebesgue}
\E\left[|u(t,x)|^2\right] =
1+(1+\vv^2)\calH(t)\;.
\end{align}
\end{corollary}

\begin{remark}\label{R2:TP-Lebesgue}
If $\rho(u)=u$ (i.e., $\lambda=1$ and $\vv=0$), then
the second moment formula \eqref{E2:SecMom-Lebesgue} recovers, in the case
$n=2$, the moment formulae of
Bertini and Cancrini
\cite[Theorem 2.6]{BertiniCancrini94Intermittence}:
\[
\E\left[|u(t,x)|^n\right] = 2 \exp\left\{\frac{n(n^2-1)}{4!\: \nu} t\right\}
\Phi\left(\sqrt{ \frac{n(n^2-1)}{12\nu} t}\right).
\]
As for the two-point correlation function, Bertini and Cancrini \cite[Corollary
2.4]{BertiniCancrini94Intermittence} gave the following integral form:
\begin{align}
\label{E2:TP-Lebesgue-BC}
\E\left[u(t,x)u\left(t,y\right)\right] =
\int_0^t \ud s \frac{|x-y|}{\sqrt{\pi \nu s^3}}
\exp\left\{-\frac{(x-y)^2}{4\nu
s}+\frac{t-s}{4\nu}\right\}\Phi\left(\sqrt{\frac{t-s}{2\nu}}\right)\:.
\end{align}
This integral can be evaluated explicitly and equals
\begin{align*}
  \exp\left(\frac{t-2 |x-y|}{4 \nu }\right)
\Erfc\left(\frac{|x-y|-t}{\sqrt{4\nu
t}}\right)\:,
\end{align*}
so their formula differs from \eqref{E2:TP-Lebesgue}. The difference is a term
$\Erf\left(|x-y|/\sqrt{4\nu t}\right)$.
By letting $x=y$ in the two-point correlation
function, both results do give the correct second moment (the
difference term is zero for $x=y$).
However, for $x\ne y$, this is not the case. For instance,
as $t$ tends to zero, the correlation
function should have a limit equal to one, while \eqref{E2:TP-Lebesgue-BC}
has limit zero.
The argument in \cite{BertiniCancrini94Intermittence} should be modified as
follows (we use the notations in
their paper): (4.6) on p. 1398 should be
\[
\E_{0}^{\beta,1}\left[
\exp\left(\frac{L_t^{\xi}(\beta)}{\sqrt{2\nu}}\right)
\right]
= \int_0^t P_\xi(\ud s)
\E_0^\beta\left[\exp\left(\frac{L_{t-s}(\beta)}{\sqrt{2\nu}}\right)\right]
+ P(T_\xi \ge t)\:.
\]
The extra term is the last term, which is
\[
P(T_\xi \ge t) = \int_t^\infty  \frac{|\xi|}{\sqrt{2\pi s^3}}
\exp\left(-\frac{\xi^2}{2s}\right) \ud s =
\Erf\left(\frac{|\xi|}{\sqrt{2t}}\right) =
\Erf\left(\frac{\left|x-x'\right|}{\sqrt{4\nu t}}\right)\;.
\]
With this term, \eqref{E2:TP-Lebesgue} is recovered.
\end{remark}

\begin{example}[Higher moments for Lebesgue initial data]
  \label{Ex2:MomLeb}
Suppose that $\mu(\ud x) = \ud x$. Clearly, $J_0(t,x)\equiv 1$.
By the above bound \eqref{E2:SecMom-Up}, we have
\begin{align*}
\E[|u(t,x)|^p] \le 
2^{p-1} +
2^{p/2-1}\left(2+\Vip^2\right)^{p/2}
\exp\left\{\frac{a_{p,\Vip}^4 \: z_p^4 \: p \: \Lip_\rho^4\:
t}{8\nu}\right\}
\left|\Phi\left(a_{p,\Vip}^2 \:\Lip_\rho^2
z_p^2\sqrt{\frac{t}{2\nu}}\right)\right|^{p/2}\:.
\end{align*}
We can replace $z_p$ by $2\sqrt{p}$, and $a_{p,\Vip}$ by
$2$. Then the upper Lyapunov
exponent of order $p$ defined in \eqref{E2:Lyapunov} is bounded by
\[
\overline{\lambda}_p \le \frac{2^5\: p^3
\Lip_\rho^4}{\nu}\:.
\]
If $\Vip=0$, we can replace $a_{p,\Vip}$ by $\sqrt{2}$ instead of $2$, which
gives a slightly better bound $\overline{\lambda}_p\le 2^3 p^3
\Lip_\rho^4/\nu$.
In particular, for the parabolic Anderson model $\rho(u) = \lambda u$, we have
\[
\overline{\lambda}_p\le 2^3 p^3 \lambda^4/\nu\;,
\]
which is consistent with Bertini and Cancrini's formulae
$\lambda_p= \frac{\lambda^4}{4! \nu} p(p^2-1)$ (see
\cite[(2.40)]{BertiniCancrini94Intermittence}).
\end{example}

\begin{corollary}[Dirac delta initial data]\label{C2:TP-Delta}
Suppose that $|\rho(u)|^2=\lambda^2(\vv^2+u^2)$ and $\mu$ is the Dirac delta
measure with a unit mass at zero. Then for all $t>0$ and $x,y\in\R$,
\begin{multline}
\label{E2:TP-Delta}
\E\left[u(t,x)u\left(t,y\right)\right] =
G_\nu(t,x)G_\nu\left(t,y\right)
-\vv^2
   \Erfc\left(\frac{|x-y|}{2 \sqrt{\nu
   t}}\right)
+\left( \frac{\lambda ^2}{4\nu}
G_{\nu/2}\left(t,\frac{x+y}{2}\right)+\vv^2\right)\\
\times
   \exp\left(\frac{\lambda ^4 t-2 \lambda ^2 |x-y|}{4 \nu }\right)
   \Erfc\left(\frac{|x-y|-\lambda ^2 t}{2
   \sqrt{\nu  t}}\right)\:.
\end{multline}
In addition, when $y=x$, we have
\begin{align}
\label{E2:SecMom-Delta}
\E\left[|u(t,x)|^2\right] =
\frac{1}{\lambda^2}\calK(t,x)+\vv^2\calH(t)\;.
\end{align}
\end{corollary}

\begin{remark}\label{R2:TP-Delta}
If $\rho(u)=u$ (i.e., $\lambda=1$ and $\vv=0$), then
the second moment formula \eqref{E2:SecMom-Delta} recovers the result
by Bertini and Cancrini
\cite[(2.27)]{BertiniCancrini94Intermittence}:
\[
\E\left[|u(t,x)|^2\right] = \frac{1}{2\pi\nu t} e^{- \frac{x^2}{\nu t}}
\left[
1+\sqrt{\frac{\pi t}{\nu}} e^{\frac{t}{4\nu}}
\Phi\left(\sqrt{\frac{t}{2\nu}}\right)
\right],
\]
which equals $\calK\left(t,x;\nu/2,1/\sqrt{4\pi \nu}\right)$.
As for the two-point correlation function, Bertini and Cancrini
\cite[Corollary 2.5]{BertiniCancrini94Intermittence} gave the following
integral form:
\begin{multline}
\label{E2:TP-Delta-BC}
 \E\left[u(t,x)u\left(t,y\right)\right] \\
= \frac{1}{2\pi\nu
t}\exp\left\{-\frac{x^2+y^2}{2\nu t}\right\}
 \int_0^1  \frac{|x-y|}{\sqrt{4\pi\nu t}}
\frac{1}{\sqrt{s^3(1-s)}}
\exp\left\{-\frac{(x-y)^2}{4\nu t}\frac{1-s}{s}\right\}
\\
\times\left(1+\sqrt{\frac{\pi t (1-s)}{\nu}} \exp\left\{\frac{t}{2\nu}
\frac{1-s}{2}\right\} \Phi\left(\sqrt{\frac{t(1-s)}{2\nu}}\right)
\right) \ud s\;.
\end{multline}
This integral can be evaluated explicitly, and is equal to
\[
 =
G_\nu(t,x) G_\nu\left(t,y\right) +
\frac{1}{4\nu}  G_{\frac{\nu}{2}}\left(t,\frac{x+y}{2}\right)
\exp\left(\frac{t-2|x-y|}{4\nu}\right)
\Erfc\left(\frac{|x-y|- t }{\sqrt{4\nu t }}\right)\:.
\]
This coincides with our result \eqref{E2:TP-Delta} for $\vv=0$ and
$\lambda=1$.
\end{remark}

\begin{example}[Higher moments for delta initial data]
\label{Ex:GrowthDeta1}
Suppose that $\mu = \delta_0$ and $\Vip=0$. Let $p\ge 2$ be an even
integer. Clearly,
$J_0(t,x)\equiv G_\nu(t,x)$. Then, by \eqref{E2:SecMom-Up}, we have that
\begin{align*}
\E\left[ |u(t,x)|^p \right]
&\le 2^{p-1} G_\nu^{p}(t,x)+2^{(p-2)/2}\left|\left(2 G_\nu^2\star
\widehat{\calK}_p\right)(t,x)\right|^{p/2}\\
&\le 2^{p-1} G_\nu^{p}(t,x)+2^{(p-2)/2}\Lip_\rho^{-p} z_p^{-p}
\left|\widehat{\calK}_p(t,x)\right|^{p/2}\\
&=
2^{p-1} G_\nu^{p}(t,x)+
2^{p-1}G_{\nu/2}^{p/2}(t,x)\left(\frac{1}{\sqrt{4\pi\nu t}}+
\frac{z_p^2\Lip_\rho^2}{\nu}\: e^{\frac{z_p^4\Lip_\rho^4
t}{\nu}}\Phi\left(z_p^2\Lip_\rho^2 \sqrt{\frac{2t}{\nu}}\right)\right)^{p/2},
\end{align*}
(the second inequality requires a proof).
Hence, for all $x\in\R$, the upper Lyapunov exponent \eqref{E2:Lyapunov} of
order $p$ is bounded by
\[
\overline{\lambda}_p \le  \frac{\Lip_\rho^4 \: z_p^4\: p}{2\nu}
\le
\frac{2^3\:  p^3 \: \Lip_\rho^4   }{\nu}\;,
\]
where the last inequality is due to the fact that $z_p\le 2\sqrt{p}$ for
all $p\ge 2$. Note that this upper bound is identical to the case of
Lebesgue initial data.
We can also calculate the exponential growth indices explicitly in this case:
\[
\lim_{t\rightarrow+\infty}\frac{1}{t} \sup_{|x|> \alpha
t}\log\E\left[|u(t,x)|^p\right]
\le -\frac{\alpha^2 p}{2\nu} + \frac{\Lip_\rho^4 \: p \: z_p^4}{2\nu}\;,\quad
\text{for all $\alpha \ge 0$}\;.
\]
Hence, the upper growth indices of order $p$ is bounded by
$\overline{\lambda}(p) \le  z_p^2 \Lip_\rho^2$. Similarly, one can derive
that $\underline{\lambda}(2)\ge \lip_\rho^2/2$.
Finally, since $\underline{\lambda}(2)\le \underline{\lambda}(p)$
for all $p\ge 2$, we have that, for all even integers $p\ge 2$,
\[
\frac{\lip_\rho^2}{2} \le \underline{\lambda}(p) \le \overline{\lambda}(p)
\le
z_p^2 \Lip_\rho^2 \:.
\]
Similar bounds are obtained for more general
initial data: see Theorem \ref{T2:Growth} below.
\end{example}

This following proposition shows that initial data cannot be extended beyond
measures.

\begin{proposition}\label{P2:D-Delta}
Suppose that the initial data is $\mu= \delta_0^{'}$, the derivative of the
Dirac delta measure at zero. Then the parabolic Anderson model $\rho(u) =
\lambda u$
($\lambda\ne 0$) does not have a random field solution in the sense of
Definition \ref{D2:Solution}.
\end{proposition}

\subsection{Exponential Growth Indices}
As an application of the above second moment formula, we partially answer
the first open problem proposed by Conus and Khoshnevisan in
\cite{ConusKhosh10Farthest}: the limits over $t$ in the definitions of these two
indices do exist when $n=2$ and the lower and upper growth indices of
order $2$ (see \eqref{E2:GrowInd-0} and \eqref{E2:GrowInd-1}) coincide.

\begin{figure}[htbp]
\centering
\includegraphics[scale=1.2]{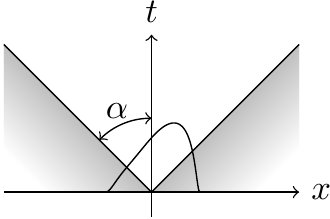}
\qquad\qquad
\includegraphics[scale=1.2]{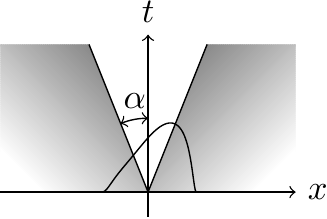}
\caption{Illustration of the exponential growth indices. The initial data,
depicted by the curve, is localized around the origin. }
\label{F2:GrowthIndices}
\end{figure}

Before stating the main result, we first give some explanation concerning the
exponential growth indices defined in \eqref{E2:GrowInd-0}
and \eqref{E2:GrowInd-1}.
When the initial data is localized, for example, when it has compact
support, we expect that the
position of high peaks of the solution will exhibit a certain wave propagation
phenomenon. As shown in Figure \ref{F2:GrowthIndices}, when $\alpha$ is
sufficiently large, it is likely that there is no high peaks outside of
the space-time cone --- the shaded region. Hence, the limit over $t$
should be negative. The largest $\alpha$ such that this limit remains negative
is then defined to be the upper growth index $\overline{\lambda}(p)$.
On the other hand, when $\alpha$ is very small, say $\alpha = 0$, then there
must be some high peaks in the shaded region so that the limit becomes
positive. Hence, the smallest $\alpha$ such that this limit is positive is
defined to be the lower growth index $\underline{\lambda}(p)$.

\begin{theorem}[Exponential growth indices] \label{T2:Growth}
The following bounds hold:
\begin{enumerate}[(1)]
 \item If $|\rho(u)|^2\le \Lip_\rho^2\left(\Vip^2+u^2\right)$ with $\Vip=0$
(which implies $\vip = \vv = 0$) and the initial data $\mu\in
\calM_{G}^{\sd}(\R)$ for some $\sd > 0$, then for all $p\ge 2$,
\[
\bar{\lambda}(p)\le
\begin{cases}
\displaystyle
\frac{\sd \nu}{2} +\frac{z_{\Ceil{p}_2}^4\Lip_\rho^4}{2\nu \sd}\:, &
\displaystyle
\text{if}\quad 0\le \sd<
\frac{z_{\Ceil{p}_2}^2\Lip_\rho^2}{\nu}\:,\cr
\displaystyle
 z_{\Ceil{p}_2}^2\: \Lip_\rho^2\;, &
\displaystyle
\text{if}\quad \sd \ge  \frac{z_{\Ceil{p}_2}^2\Lip_\rho^2}{\nu}\:,
\end{cases}
\]
where $z_m$, $m\in \bbN$, $m\ge 2$, are the universal constants in the
Burkholder-Davis-Gundy inequality.
In particular, for $p=2$,
\begin{align}
\bar{\lambda}(2)\le
\begin{cases}
\displaystyle
\frac{\sd\nu}{2} +\frac{\Lip_\rho^4}{8\nu \sd}\;, &
\displaystyle
\text{if}\quad 0\le \sd< \frac{\Lip_\rho^2}{2\nu}\;,\cr
\displaystyle
 \frac{1}{2}\Lip_\rho^2\;, &
\displaystyle
\text{if}\quad \sd \ge  \frac{\Lip_\rho^2}{2\nu}\:.
\end{cases}
\label{E2:Growth-L2}
\end{align}
\item If $|\rho(u)|^2\ge \lip_\rho^2\left(\vip^2+u^2\right)$ with $\vip=0$, then
\[
\underline{\lambda}(p)\ge \frac{\lip_\rho^2}{2},\qquad \text{for all $\mu\in
\calM_{+}(\R)$, $\mu\ne 0$ and all $p\ge 2$}\:;
\]
otherwise, if $\vip\ne 0$, then
\[
\underline{\lambda}(p) =
\overline{\lambda}(p)=+\infty,\qquad \text{for all $\mu\in
\calM_{+}(\R)$ and $p\ge 2$}\:;
\]
\item In particular, for the quasi-linear case
$|\rho(u)|^2=\lambda^2\left(\vv^2+u^2\right)$
with $\lambda\ne 0$, if $\vv=0$ and $\sd \ge \frac{\lambda^2}{2\nu}$, then
\[
\underline{\lambda}(2)=\bar{\lambda}(2)=\lambda^2/2,\qquad \text{for all
$\mu\in \calM_{G,+}^{\sd}(\R)$, $\mu\ne 0$}\:;
\]
otherwise, if $\vv\ne 0$, then
\[
\underline{\lambda}(p) =
\overline{\lambda}(p)=+\infty,\qquad \text{for all $\mu\in
\calM_{+}(\R)$ and $p\ge 2$}\;.
\]
\end{enumerate}
\end{theorem}

This theorem generalizes the results by Conus
and Khoshnevisan \cite{ConusKhosh10Farthest} in several aspects:
(i) more general initial data are allowed; (ii) both non trivial upper bound and
lower bounds are given (compare with Theorem 1.1 \cite{ConusKhosh10Farthest})
for the Laplace operator case;  (iii) for the parabolic Anderson model,
the exact transition is proved (see Theorem 1.3 and the first open problem in
\cite{ConusKhosh10Farthest}) for $n=2$ and the Laplace operator case; (iv)
our discussions above cover the case $\rho(0)\ne 0$.

\begin{example}[Dirac delta initial data]
Suppose that $\Vip=\vip=0$. Clearly, $\delta_0 \in \calM_{G,+}^{\sd}(\R)$ for
all $\sd\ge 0$. Hence, the above
theorem implies that for all even integers $k\ge 2$,
\[
\frac{\lip_\rho^2}{2} \le \underline{\lambda}(k) \le \overline{\lambda}(k) \le
z_k^2 \Lip_\rho^2\;.
\]
This recovers the previous calculation in Example \ref{Ex:GrowthDeta1}.
\end{example}


\begin{proposition}\label{P2:Example-Exponents}
Consider the parabolic Anderson model $\rho(u) =\lambda u$, $\lambda\ne 0$ with
the initial data $\mu(\ud x) = e^{-\sd|x|}\ud x$ ($\sd>0$). Then we have
\[
\underline{\lambda}(2)=\overline{\lambda}(2) =
\begin{cases}
\displaystyle
\frac{\sd\nu}{2} +
\frac{\lambda^4}{8 \sd \nu}& \text{if $\;\;\displaystyle 0<\sd\le
\frac{\lambda^2}{2\: \nu}$}\;,\cr\cr
\displaystyle
\frac{\lambda^2}{2}
& \text{if $\;\;\displaystyle \sd\ge
\frac{\lambda^2}{2\: \nu}$}\;.
\end{cases}
\]
\end{proposition}

This proposition shows that for all $\sd\in
\left]0,+\infty\right]$, the exact phase transition occurs,
and hence our upper bounds \eqref{E2:Growth-L2} in Theorem \ref{T2:Growth} for
the upper growth index
$\overline{\lambda}(2)$ are sharp.

\subsection{Sample Path Regularity}

\begin{theorem}\label{T2:Holder}
Suppose that $\rho$ is Lipschitz continuous. Then the solution
$u(t,x)=J_0(t,x)+I(t,x)$ to \eqref{E2:Heat} has the following sample path
regularity:
\begin{enumerate}[(1)]
 \item If the initial data $\mu$ is an $\alpha$-H\"older continuous function
($\alpha\in \;]0,1]$) over $\R$ satisfying \eqref{E2:J0finite}, then
\[
J_0\in
C_{\frac{1}{2},\alpha}
\left(\R_+\times\R\right)
\;\cup\;
C_{\frac{1}{2},1}\left(\R_+^*\times\R\right)\;,
\quad\text{and}\quad
I\in C_{\frac{1}{4}-,\frac{1}{2}-}\left(\R_+\times\R\right),\quad \text{a.s.}
\]
Therefore,
\[
u=J_0+I\in
C_{\frac{1}{4}-,\left(\frac{1}{2}-\right)\wedge
\alpha}\left(\R_+\times\R\right)
\;\cup\;
C_{\frac{1}{4}-,\frac{1}{2}-}\left(\R_+^*\times\R\right)
,\quad \text{a.s.}
\]
\item If the initial data $\mu$ is a continuous function
satisfying \eqref{E2:J0finite}, then
\[
J_0\in
C_{\frac{1}{2},1}\left(\R_+^*\times\R\right)\;,
\quad\text{and} \quad
I \in C_{\frac{1}{4}-,\frac{1}{2}-}\left(\R_+\times\R\right),\quad \text{a.s.}
\]
Therefore,
\[
u=J_0+I\in
C_{\frac{1}{4}-,\frac{1}{2}-}\left(\R_+^*\times\R\right)
,\quad \text{a.s.}
\]
\item If the initial data $\mu$ is a signed Borel measure satisfying
\eqref{E2:J0finite}, then
\[
J_0\in
C_{\frac{1}{2},1}\left(\R_+^*\times\R\right)\;,
\quad\text{and}\quad
I\in C_{\frac{1}{4}-,\frac{1}{2}-}\left(\R_+^*\times\R\right),\quad \text{a.s.}
\]
Therefore,
\[
u=J_0+I\in
C_{\frac{1}{4}-,\frac{1}{2}-}\left(\R_+^*\times\R\right)
,\quad \text{a.s.}
\]
\end{enumerate}
\end{theorem}

\begin{example}[Dirac delta initial data]\label{E2:Holder-Delta}
Suppose $\rho(u)= \lambda u$ with $\lambda \ne 0$.
If $\mu=\delta_0$, then neither $J_0(0,x)$ nor $\lim_{t\rightarrow 0_+}
\Norm{I(t,x)}_2$ is continuous in $x$.
For $J_0(0,x) = \delta_0(x)$, this is clear.
As for $\lim_{t\rightarrow 0_+} \Norm{I(t,x)}_2$, by Corollary \ref{C2:TP-Delta}
(with $\vv=0$), we have
\[
\Norm{I(t,x)}_2^2=  \frac{1}{\lambda^2}\calK(t,x) -G_\nu^2(t,x)
= \frac{\lambda^2}{2\nu}e^{\frac{\lambda^4
t}{4\nu}}\Phi\left(\lambda^2\sqrt{\frac{t}{2\nu}}\right) G_{\nu/2}(t,x)\;.
\]
Therefore,
\[
\lim_{t\rightarrow 0_+}
\Norm{I(t,x)}_2^2 =
\begin{cases}
 0& \text{if $x\ne 0$}\;,\cr
+\infty & \text{if $x=0$}\;.
\end{cases}
\]
\end{example}

\bibliographystyle{abbrv}

\end{document}